\documentclass[12pt,letterpaper]{article}

\newcommand{\required}[1]{\section*{\hfil \sharp1\hfil}}

\usepackage{verbatim,url}
\usepackage{graphicx}

\usepackage{amssymb,amsmath,amsfonts,amsthm}
\usepackage{hyperref}
\usepackage{upref}
\newcommand{\Beq}{\begin{equation}}
\newcommand{\Eeq}{\end{equation}}
\newcommand{\beq}{\begin{equation*}}
\newcommand{\eeq}{\end{equation*}}
\newcommand{\bal}{\begin{align}}
\newcommand{\eal}{\end{align}}

\newtheorem{theorem}{Theorem}
\newtheorem{Lemma}{Lemma}

\newtheorem{definition}{Definition}

\theoremstyle{definition}

\newtheorem{remark}{Remark}

\newcommand{\ac}[1]{\begin{quotation}\textbf{Anuj's comment:\
		}{\textit{\sharp1}}\end{quotation}}
\newcommand{\rc}[1]{\begin{quotation}\textbf{Rohit's comment:\
		}{\textit{\sharp1}}\end{quotation}}
\allowdisplaybreaks

\usepackage{xcolor}

\title{Adaptive estimation of a function from its Exponential Radon Transform in presence of noise}
\date{}
\author{Anuj Abhishek$^{a}$ and Sakshi Arya$^{b}$\\[1ex]
	\footnotesize $^{a}$ Dept. of Mathematics, Drexel University, Philadelphia, PA, USA (\emph{anuj.abhishek@drexel.edu}) \hfill\mbox{}\\[-0.5ex]
	\footnotesize $^{b}$ Dept. of Public Health Sciences, University of Chicago, IL, USA (\emph{sakshi@uchicago.edu}) \hfill\mbox{}
}

\begin{document}
	\maketitle
	\begin{abstract}
		In this article we propose a locally adaptive strategy for estimating a function from its Exponential Radon Transform (ERT) data, without prior knowledge of the smoothness of  functions that are to be estimated. We build a non-parametric kernel type estimator and show that for a class of functions comprising a wide Sobolev regularity scale, our proposed strategy follows the minimax optimal rate up to a $\log{n}$ factor. We also show that there does not exist an optimal adaptive estimator on the Sobolev scale when the pointwise risk is used and in fact the rate achieved by the proposed estimator is the adaptive rate of convergence. 
	\end{abstract}
	\section{Introduction}
\label{Intro} Single Photon Emission Computed Tomography (SPECT) imaging is a valuable diagnostic tool that is frequently used to detect the presence of tumors inside a patient's body. The idea behind SPECT imaging can be described very briefly in the following manner: A small amount of radioactive tracer attached to some nutrient is injected in the patient's body. After a brief interlude (ranging from a few minutes to a few hours), a SPECT scanner is used to measure the radioactive emissions from the body in a range of directions by moving the scanner around the body. Along each line, the data represents the intensity of emissions from a point along that line. This data can be mathematically interpreted as an attenuated Radon transform. From the attenuated Radon transform data, one then tries to image the inside of the patient's body to locate the presence of tumors. If one makes the simplifying assumption that the attenuation is constant, then the attenuated Radon transform reduces to the case of what is known as the exponential Radon transform. We point the interested reader to \cite{Kuchment_book} and \cite {natt_wubb_book} for a more detailed overview.
\par \noindent In the setting of the current article, our focus of investigation is the estimation of a function from its stochastic (i.e. noisy) exponential Radon transform (ERT) data. In fact, the ERT of a compactly supported function $f(x)$ in $\mathbb{R}^2$ is given by:
\begin{align}
T_{\mu}f(\theta,s)=\int\limits_{x\cdot\theta=s}e^{\mu x\cdot\theta^{\perp}}f(x)dx.
\end{align}
Here $s\in \mathbb{R}$, $\theta\in \mathrm{S}^1$ where $\mathrm{S}^1$ is the unit circle in $\mathbb{R}^2$, $\mu$ is a constant and $\theta^{\perp}$ denotes a unit vector perpendicular to $\theta$. Recall that lines in $\mathbb{R}^2$ can be parameterized as $L(\theta,s)=\{x:x\cdot\theta=s\}$. Thus, just as the classical Radon transform, ERT takes a function defined on a plane and maps it to a function defined over the set of lines parameterized by $(\theta,s)$. Indeed, the attenuated Radon transform (and thus, the ERT) is itself an example of generalized Radon transforms that were studied in \cite{Quinto_80, Quinto_83}.  Inversion methods for the exponential Radon transform (in a non-noisy setting) are known from \cite{Natt_79} and \cite{TM_80}, see also \cite{hazou_solmon} for filtered backprojection (FBP) type formulas. 
\par \noindent Classical Radon transform has also been extensively studied in the stochastic framework. The problem of positron emission tomography (PET) in presence of noise was studied in \cite{JS_90}. In \cite{Tsybakov_1991,Tsybakov_92,minimax_book} the authors show that the kernel type non-parametric estimators (which are closely linked to FBP inversion methods) attain optimal minimax convergence rate. In \cite{AA_opt_ERT}, the author extended the results that were known for Radon transform from \cite{Tsybakov_1991,Tsybakov_92} to the setting of stochastic ERT. In \cite{Cavalier_98} Cavalier obtained results on efficient estimation of density in the non-parametric setting for stochastic PET problem. In addition to the non-parametric kernel type estimators, Bayesian estimators for the stochastic problem of X-ray tomography have been studied by several authors, see e.g. \cite{Siltanen2003,Lassas_09} and references therein. More recently, authors in \cite{Monard_19} have obtained results on efficient Bayesian inference for the attenuated X-ray transform on a Riemannian manifold. In the tomography results described so far, in the stochastic framework one usually assumes that the the smoothness (e.g. Sobolev regularity) of the function to be estimated is known. An interesting problem is to devise adaptive estimation procedures that can be applied for the estimation of a function without a priori knowledge of its smoothness. In \cite{Cavalier01}, the problem of estimation of bounded functions from its noisy Radon transform data was solved. The (locally adaptive) estimation procedure described in \cite{Cavalier01} was based on the method proposed in \cite{Lepski_97}. Adaptive estimation of functions has a rich and varied history and was given a major impetus by publication of a series of articles \cite{Lepski90,Lepski91,Lepski92}. Spatially adaptive estimation procedures were considered in \cite{Donoho94,Lepski95}. The problem of optimal pointwise adaptation was considered in \cite{Lepski97a, Tsybakov_98} for H\"older and Sobolev classes respectively. In fact, an ubiquitous feature of pointwise adaptive estimation over Sobolev classes seems to be a certain logarithmic loss of efficiency when {compared to the optimal minimax rate of convergence if} one assumes knowledge about the Sobolev regularity of the function, see also \cite{Butucea_00, Butucea_01} for related results. Such a loss of efficiency for the estimation of a function from its Radon transform data was also conjectured (on the Sobolev scale) by the author in \cite{Cavalier01}. Our results in the current article establish that it indeed is the case in the even more general set-up of estimation of function from its exponential Radon transform. We remark here that the adaptive estimation of function from its ERT data falls in the category of statistical inverse problems which are somewhat different in flavor from the problems in adaptive estimation in the non-parametric regression setting or probability density estimation from direct observations. Finally, we mention below a partial list of several important results in the theory of adaptive estimation for inverse problems such as deconvolution and change point estimation which have pushed the boundaries of this area of research far and wide, see e.g. \cite{Tsybakov07,Tsybakov07a,Tsybakov03,Tsybakov02,gold99,Tsybakov08a,Tsybakov08,Lepski17, Lepski19}.

The organization of this article is as follows: in section 2, we describe the mathematical set-up of the problem and recall some relevant definitions. In section 3, we apply the adaptive strategy of \cite{Cavalier01} to the problem of stochastic ERT and at first recover similar results for the ERT case as was proved in the Radon transform setting by the author in \cite{Cavalier01}. The proofs for the first three theorems follow from similar techniques as in \cite{Cavalier01} except for the modifications needed to adapt the proofs to the ERT problem. We present these proofs for the sake of completeness. Theorem \ref{Th4} in the article establishes a `no-optimality' result for the adaptive estimation of a function from its ERT data and shows that among all the adaptive strategies, the strategy as used in the current article is the `best' (see definition \ref{adap}). For the proof of this theorem, we have adapted the method used in \cite{Butucea_00,Butucea_01} in the direct problem of density estimation over Sobolev classes in the setting of our particular inverse problem. Finally, the appendix has a proof of an auxiliary lemma.
\section{Problem set-up and definitions}\label{section2}
\label{setup}
\noindent In this section we will describe the mathematical framework for the problem and recall some standard definitions from the literature. 

Let {$B_1=\{x:\vert\vert x\vert\vert \leq 1\}$} be the unit ball in $\mathbb{R}^2$. Let $f(x):\mathbb{R}^2\to \mathbb{R}$ be a function such that it is supported in {$B_1$}, is continuous (a.e.) and $\lvert f(x)\rvert \leq L$ for some $L>0$. We will denote the class of such functions by $B(L)$.\\
\begin{definition}
	Let $\mathrm{S}^1$ denote the unit circle in $\mathbb{R}^2$ and $Z= \mathrm{S}^1\times [-1,1]$ be the cylinder whose points are given by $(\theta,s)$ where $s\in [-1,1]$ and $\theta \in \mathrm{S}^1$. By $\theta^{\perp}$, we will denote a unit vector perpendicular to $\theta$. The exponential Radon transform of $f\in B(L)$ is defined as the following function on $Z$:
	$$T_{\mu}f(\theta,s)=\int_{x\cdot\theta=s}e^{{\mu}x\cdot\theta^{\perp}}f(x)dx,$$ where $\mu$ is a fixed constant. It is clear that if $\mu=0$, then the exponential Radon transform reduces to the case of the classical Radon transform. 
\end{definition}
\begin{definition}
	Let $g(\theta,s)$ be a compactly supported function on $Z$. The associated dual transform is then defined as: $$T_\mu^{\sharp}g(x)=\int_{\mathrm{S}^1}e^{\mu x\cdot\theta^\perp}g(\theta,x\cdot \theta)d\theta.$$ For $\mu=0$, this reduces to the backprojection operator for the classical Radon transform.
	
\end{definition}
\noindent  Let $\{(\theta_i,s_i)\}_{i=1}^{i=n}$ be $n$ random points on the observation space $Z$ and let the observations be of the form:
\begin{equation} \label{obs_model}
Y_i=T_{\mu}f(\theta_i,s_i)+\epsilon_i.
\end{equation} 
\par \noindent We assume that the points $(\theta_i,s_i)$ are independent and identically distributed (i.i.d.) on $Z$ and are distributed uniformly. The noise, $\epsilon_i$, are  i.i.d. Gaussian random variables with zero mean and some finite positive variance $\sigma^2$. 
Then the stochastic inverse problem for exponential Radon transform is to estimate the function $f(x)$ based on the observations $Y_i$ for $i=\{1,2,\dots,n\}$. Let us denote by $\hat{f}_n(x)$ any estimator of $f(x)$ based on the observed data. 

Now we recall some definitions that will be frequently used in this article. In this article, the semi-norm $d$ (e.g. in Definition \ref{risk_fn}) will be given by: $d(f,g)=\vert f(x_0)-g(x_0)\vert$ where $x_0$ is an arbitrary fixed point in  {$B_1$}. From here on, $E_f(\cdot)$ and $Var_f(\cdot)$ will be used to denote the expectation and variance with respect to the joint distribution of random variables $(s_i,\theta_i,Y_i)$, $i=\{1,\dots,n\}$ satisfying the model given by (\ref{obs_model}).
\begin{definition}[\cite {JS_90,Tsybakov_book}]\label{risk_fn}
	The risk function of an estimator $\hat{f}_n(x)$ is defined as:
	$$ \mathcal{R}(\hat{f}_n,f)=E_f(d^2(\hat{f}_n,f)).$$ 
	When the $d$ is as above, the risk is also referred to as the Mean Squared Error (MSE).
\end{definition}
\noindent The overall measure of risk associated to any estimation procedure is given by the \textit{minimax risk}.
\begin{definition}\cite[Page 78]{Tsybakov_book}\label{minimax_risk}
	Let $f(x)$  belong to some non-parametric class of functions $\mathcal{F}$. The maximum risk of an estimator $\hat{f}_n$ is defined as:
	$$r(\hat{f}_n)=\sup_{f\in\mathcal{F}} \mathcal{R}(\hat{f}_n,f). $$ Finally, the minimax risk on $\mathcal{F}$ is defined as:
	$$r_n(\mathcal{F})=\inf_{\hat{f}_n} \sup_{f\in\mathcal{F}}\mathcal{R}(\hat{f}_n,f),$$ where the infimum is taken over the set of all possible estimators $\hat{f}_n$ of $f$. Clearly, $$r_n(\mathcal{F})\leq r(\hat{f}_n).$$
\end{definition}
\noindent In the next definition, we recall the concept of minimax optimality.
\begin{definition} \label{Def5} \cite[Page 78]{Tsybakov_book}\label{optimality} Let $\{\Psi_n^2\}_{n=1}^{\infty}$ be a positive sequence converging to zero. An estimator $\hat{f}_n^{*}$ is said to be minimax optimal if there exist finite positive constants $C_1$ and $C_2$ such that,
	$$ C_1\Psi_n^2\leq r_n(\mathcal{F})\leq r(\hat{f}_n^{*})\leq C_2\Psi_n^2. $$
	Furthermore, $\Psi_n^2$ is said to be the optimal rate of convergence.
	
\end{definition}
\par \noindent 

A function $f$ belongs to the Sobolev ball $H(\beta,C)$, if $$\int_{\mathbb{R}^2}(1+\vert\vert \xi \vert\vert^2)^\beta {\vert \widetilde{f}(\xi)\vert^2} d\xi\leq C,$$ where $\widetilde{f}(\xi)$ denotes the Fourier transform of $f$. We will assume $\beta>1$ and wherever understood, we will write $H(\beta)$ for $H(\beta,C)$. In \cite{AA_opt_ERT}, we showed that the minimax optimal rate of convergence under the MSE risk in the estimation of a function on $\mathbb{R}^2$ from its stochastic ERT is given by, {$\phi_{n,\beta}=\mathcal{O}(n^{-({\beta-1})/({2\beta+1})})$}. This follows as a consequence of Theorems 3 and 5 in \cite{AA_opt_ERT}. 
Note that such an optimal minimax rate can be achieved by an estimator if one knows the smoothness of the function that is to be estimated (in particular that it belongs to $H(\beta, C)$), however in practice the smoothness is unknown.
\par \noindent In this article, our goal is to build locally adaptive data driven estimators that do not assume prior knowledge about the smoothness of the functions that are to be estimated. We will test the accuracy of such estimators by looking at their performance over a class of functions encompassing a wide scale of Sobolev regularity. Let us make these ideas precise: assume now that we only know that the function to be estimated belongs to $H(\beta)\cap B(L)$ where $\beta$ lies in a discrete set $B_n$ given by $B_n=\{\beta_1<\dots<\beta_{N_n}\}$ such that $\beta_1>1$ is fixed and $\lim_{n\to\infty}B_{N_n}\to \infty$.  The adaptive rate of convergence (ARC) on a scale of classes $H(\beta)\cap B(L)$, $\beta\in B_n$ is defined as:
\begin{definition} \cite[Definition 3]{Tsybakov_98}\label{adap}
	A sequence $\psi_{n,\beta}$ is said to be an ARC if:
	
	(a) There exists a rate adaptive estimator $f^*(x)$ independent of the smoothness scale $\beta$ such that
	\begin{align}
	\limsup_{n\to \infty}\sup_{\beta \in B_n}\sup_{f\in H(\beta)\cap B(L)}(\psi_{n,\beta})^{-2}E_{f}(f^*(x)-f(x))^2<\infty. \label{adaptive_conditionA}
	\end{align}
	
	(b) If there exists another sequence $\gamma_{n,\beta}$ and another adaptive estimator $f^{**}(x)$ satisfying:
	$$\limsup_{n\to \infty}\sup_{\beta \in B_n}\sup_{f\in H(\beta)\cap B(L)}(\gamma_{n,\beta})^{-2}E_{f}(f^{**}(x)-f(x))^2<\infty.$$
	and a $\beta^{\prime }$ such that $$\frac{\gamma_{n,\beta^{\prime}}}{\psi_{n,\beta^{\prime}}}\underset{n \to \infty}{\to} 0,$$ then there exists a $\beta^{\prime\prime}$ such that $$\frac{\gamma_{n,\beta^{\prime}}}{\psi_{n,\beta^{\prime}}}\frac{\gamma_{n,\beta^{\prime\prime}}}{\psi_{n,\beta^{\prime\prime}}}\underset{n \to \infty}{\to} \infty.$$
\end{definition}
{In other words, if some rate other than $\psi_{n,\beta}$ satisfies a condition similar to \eqref{adaptive_conditionA} and if this rate is faster for some smoothness parameter $\beta^\prime$, then there has to be some other smoothness parameter, $\beta^{\prime \prime}$, where the loss is infinitely greater for large sample sizes $n$.}
\begin{remark}\label{orc}
	An adaptive estimator is said to be optimally rate adaptive if it achieves minimax optimality for every $\beta \in B_n$, see \cite  [equation 2.6]{Tsybakov_98}. If there exists an estimator that is optimally rate adaptive then it also achieves the adaptive rate of convergence.
\end{remark}
Next, we discuss the procedure for building an adaptive strategy and present results that establish the adaptive rate of convergence of the proposed strategy.

\section{Adaptive strategy}
We begin by recalling some results from \cite{AA_opt_ERT}. Let $$K_{\delta}(s)=\frac{1}{\pi}\int_{\vert\mu\vert}^{\sqrt{(1/\delta^2)+\mu^2}}r \cos(sr)dr.$$ These kind of functions have been used in the context of filtered backprojection formulas for Radon transforms, see e.g. \cite[Page 237]{minimax_book}, \cite [Page 109]{Natterer_textbook} and the quantity $\delta$ is referred to as the bandwidth of the filter.  Let $\star$ represent the operation of convolution of functions. Furthermore, whenever the convolution of two functions $f$ and $g$ defined on the cylinder $Z= \mathrm{S}^1\times \mathbb{R}$  is considered, the convolution will be understood to be taken with respect to their second variable, i.e. $f\star g(\theta,s)=\int_{\mathbb{R}}f(\theta,s-t)g(\theta,t)dt.$ We consider the estimator:
\begin{align}
\bar{f}_{\delta_n}(x)=\frac{1}{n}\sum_{i=1}^n e^{-\mu x\cdot \theta_i^{\perp}}K_{\delta_n}(\langle x\cdot \theta_i \rangle -s_i)Y_i,
\end{align}
where $Y_i$ is the observed data as in equation \ref{obs_model}. For this estimator we evaluated the bias $B(\bar{f}_{\delta_n})\leq{\tilde{c}}{\delta_n^{\beta-1}}$ where $\tilde{c}$ is a constant as well as its variance {$Var_f(\bar{f}_{\delta_n})\leq c^{*}(n\delta_n)^{-3} =v^2(\delta_n) $} in \cite{AA_opt_ERT}. By balancing the bias and the variance terms it was shown that if $\delta_n=c_0\cdot[n^{{-1}/{(2\beta+1)}}]$ where $c_0$ is a constant, then the estimator is minimax optimal. Notice that the choice of such an optimal bandwidth depends upon the smoothness $\beta$ of the function to be estimated. We will now describe an adaptive bandwidth selection procedure that can be used when the smoothness of function to be estimated is not known. This bandwidth selection procedure was proposed by Cavalier \cite{Cavalier01} in the context of adaptive estimation of a function from its stochastic Radon transform data and is based on the method described in the article by Lepski \textit{et al.} \cite{Lepski_97}.
\par
\noindent We assume that the locally adaptive bandwidths {$\bar{\delta}(x)$} are chosen from a geometrical grid $\Delta_n$ given by:
\begin{align}
\Delta_n =\{\delta \in [\delta_n^{-},1]:\quad  \delta=a^{-j}, \quad j=0,1,2,\dots \}, \label{Delta_n_def}
\end{align}
where $a \geq 2$ {and $\frac{a \log{n}}{n} \leq 1, \delta_n^{-} = \frac{\log{n}}{n}$}. Let $f_{\delta} {(x)}= E_f[\bar{f}_{\delta}(x)]$. Furthermore, similar to the proof of \cite[Lemma 1]{AA_opt_ERT} and \cite [equation 12]{Cavalier01}, it can be shown that for some constant $c^{**}$, \begin{align*}
\text{Var}_f (\bar{f}_{\delta}(x)-\bar{f}_{\eta}(x))\leq \frac{c^{**}}{n}\int (K_\delta(u)-K_{\eta}(u))^2 du:=v^2(\delta,\eta).
\end{align*}  For $\delta>\eta$, we define,
\begin{align}
\psi(\delta,\eta)=v(\delta)\lambda(\delta)+v(\delta,\eta)\lambda(\eta),
\end{align}
where $\lambda(\delta)=\max (1,\sqrt{D_2 \log \frac{1}{\delta}})$ and {$v^2(\delta) = c^{*}(n\delta)^{-3} \geq Var_f(\bar{f}_{\delta})$}. Here, $D_2$ is a real number which can be chosen as desired (we will make the choice more precise later). The data driven `adaptive bandwidth' $\bar{\delta}(x)$ will be given by the following relation:
$$\bar{\delta}(x)=\max\{\delta\in \Delta_n: \lvert \bar{f}_{\delta }(x)-\bar{f}_{\eta}(x) \rvert \leq \psi(\delta,\eta)\quad \forall \eta\leq \delta, \eta\in \Delta_n \}.$$
Correspondingly, the adaptive estimator is given by: \begin{align}\label{ad_es}
f^*(x)=\bar{f}_{\bar{\delta}}(x)=\frac{1}{n}\sum_{i=1}^{n}e^{-\mu x\cdot \theta_i^{\perp}}K_{\bar{\delta}(x)} (\langle \theta_i,x \rangle -s_i)Y_i.
\end{align} We remark here that the definition for $\bar{\delta}(x)$ is defined locally at every point $x$ and it does not assume the a priori knowledge about the function $f(x)$, in particular its smoothness. Next, we define a `locally deterministic bandwidth' whose definition involves the unknown function $f(x)$ itself: 
\begin{align} \label{det_bw}
\delta_n=\delta_n(x,f)=\max\{\delta\in \Delta_n:\lvert f_{\eta}(x)-f(x)\rvert \leq \frac{v(\delta)\lambda (\delta)}{2} \quad \forall \eta\in \Delta_n, \eta\leq \delta \}.
\end{align}
Finally, following \cite{Cavalier01,Lepski_97} we define the adaptive convergence rate $r_n(x,f)$:
\begin{align}
r_n(x,f) =\inf_{\delta \in [\delta_n^-,1]}\bigg\{ \sup_{0\leq\eta\leq\delta}( f_{\eta}(x)-f(x))^2+c^*\delta^{-3}\log n/n\bigg\}.
\end{align}
Now we are ready to state our first theorem which essentially says that the `estimator' formed with `locally deterministic bandwidth' $\delta_n$ has its risk bounded by $r_n(x,f)$ up to a constant factor. However, we also note that the $\delta_n$ can be found only if one knows the function $f(x)$ in the first place. In this sense, $\delta_n$ can be thought of as an ideal bandwidth for the adaptive estimation procedure {and $\bar{f}_{\delta_n}(x)$ as an oracle}. 
\begin{theorem}
	For any $f\in B(L)$ we have as $n\to \infty$:
	\begin{align}
	E_f[(\bar{f}_{\delta_n}(x)-f(x))^2]\leq \frac{5}{4} v^2(\delta_n)\lambda^2(\delta_n)\leq C(a) r_n(x,f),
	\end{align}
	where $C(a)$ is a constant depending on $a$, {where $a$ is the same as in \eqref{Delta_n_def}}.
\end{theorem}
\begin{proof}
	
	\noindent First of all, we address an auxiliary point. In order for the definition of $ \delta_n$ to make sense, we will first need to show that the set over which the maximum is taken in this definition is non-empty. 
	Recall that,
	$f_{\delta_n^-}(x)=\boldsymbol{\delta}^{\frac{1}{\delta_n^-}}\star f$ (see eg. \cite{AA_opt_ERT}) where, $$\boldsymbol{\delta}^{\frac{1}{\delta_n^-}}=\int_{\lvert \xi \rvert \leq \frac{1}{\delta_n^-}}e^{-i\xi\cdot x}d\xi= \int I_{\frac{1}{\delta_n^-}}(\xi)e^{-i\xi\cdot x}d\xi.$$ Thus,
	\begin{align}\label{bound1}
	\lvert (f_{\delta_n^{-}}&-f)(x)\rvert \leq 2\lvert f_{\delta_n^{-}}(x)\rvert ^2+ 2\lvert f(x)\rvert ^2 =2\lvert \boldsymbol{\delta}^{\frac{1}{\delta_n^-}}\star f(x)\rvert ^2+ 2\lvert f(x)\rvert ^2 \nonumber \\
	&\leq 2\bigg(\frac{1}{4\pi^2}\bigg)\bigg(\int_{\mathbb{R}^2}\lvert\hat{f}(\xi)I_{\frac{1}{\delta_n^-}}(\xi)\rvert d\xi\bigg)^2 +2 \lvert f(x)\rvert ^2\nonumber \\
	&\leq  2\bigg(\frac{1}{4\pi^2}\bigg)\bigg(\int_{\mathbb{R}^2}\lvert\hat{f}(\xi)\rvert ^2 d\xi \int_{\mathbb{R}^2}\lvert I_{\frac{1}{\delta_n^-}}(\xi)\rvert d\xi\bigg) +2 \lvert f(x)\rvert ^2 \quad (\text{H\"{o}lder's inequality}) \nonumber\\
	&\leq  2\bigg(\frac{1}{4\pi^2}\bigg) L^2 \pi \cdot \pi \bigg(\frac{1}{\delta_n^-}\bigg)^2+2L^2 \leq \frac{5 L^2}{{2 \delta_n^-}^2}.
	\end{align}
	On the other hand, $\frac{1}{4}v^2(\delta_n^-) \lambda^2(\delta_n^-)\geq \frac{1}{4}c^*\frac{(\delta_n^-)^{-3}}{n} D_2 \log (\frac{1}{\delta_n^-})$. We also have, $\log (\frac{1}{\delta_n^-})=\log(\frac{n}{\log n})\geq \frac{\log n}{2}=n\delta_n^-/2$. Thus,  \begin{align}
	\frac{1}{4}v^2(\delta_n^-) \lambda^2(\delta_n^-)\geq \frac{c^* D_2}{8 (\delta_n^-)^2}. \nonumber
	\end{align}
	Note that if $D_2\geq 20 L^2/c^*$, then $\frac{c^* D_2}{8 (\delta_n^-)^2}\geq (5 L^2)/(2 (\delta_n^-)^2).$ This along with (\ref{bound1}) gives, $$\lvert f_{\delta_n^-}(x)-f(x)\rvert \leq \frac{1}{4}v^2(\delta_n^-)\lambda^2(\delta_n^-),$$ which in turn shows that the set over which the maximum is taken in {definition \ref{det_bw}} is non-empty. Now coming back to the proof the theorem,
	\begin{align}
	E_f(\bar{f}_{\delta_n}(x)-f(x))^2&=(f_{\delta_n}(x)-f(x))^2+\text{var}_f \bar{f}_{\delta_n}(x) \nonumber\\
	&\leq \frac{1}{4}v^2(\delta_n)\lambda^2(\delta_n)+v^2(\delta_n)\nonumber\\
	&\leq \frac{5}{4}v^2(\delta_n)\lambda^2(\delta_n) \quad (\text{using }\lambda^2(\delta_n)\geq 1). \nonumber
	\end{align} 
	Let the infimum in the above definition of $r_n(x,f)$ be obtained for $\delta=\delta_0$. We now have two cases:\\
	\textbf{Case 1: }If $\delta_0<a\delta_n$,
	\begin{align}
	r_n(x,f)\geq c^*\delta_0^{-3}\frac{\log n}{n}\geq c^* a^3 \delta_n^{-3}\frac{\log n}{n}. \nonumber
	\end{align}
	From the definitions of $v^2(\delta_n)$ and $\lambda^2(\delta_n)$, we know that, $$\frac{5}{4}v^2(\delta_n)\lambda^2(\delta_n)=\max\bigg\{\frac{5}{4}\frac{c^*\delta_n^{-3}}{n} , \frac{5}{4}\frac{c^*\delta_n^{-3}D_2}{n}\log(1/\delta_n) \bigg\}.$$ Since $\delta_n^-\leq \delta_n$, then for $n\geq 3$ we have the following sequence of inequalities: $$\log (1/\delta_n)\leq \log (1/\delta_n^-) = \log (n/\log n) \leq \log n.$$ Thus if we choose $C_1=C_1(a,D_2)>\max\{\frac{5}{4 a^3}, \frac {5 D_2}{4 a^3} \}$, we have:
	\begin{align}
	E_f[(\bar{f}_{\delta_n}(x)-f(x))^2]\leq \frac{5}{4}v^2(\delta_n)\lambda^2(\delta_n)\leq C_1 r_n(x,f), \nonumber
	\end{align}
	for $n$ large enough.\\
	\textbf{Case 2: } If $\delta_0\geq a\delta_n$, then from the definition of $\delta_n$ (see (\ref{det_bw})),
	$$ \sup_{0\leq \eta\leq \delta_0}(f_{\eta}(x)-f(x))^2 \geq \sup_{0\leq \eta\leq a \delta_n}(f_{\eta}(x)-f(x))^2\geq \frac{v^2(a\delta_n)\lambda^2(a\delta_n)}{4}.$$
	So for some $C_2=C_2(a)$ large enough, we have:\\ $$E_f[(\bar{f}_{\delta_n}(x)-f(x))^2]\leq \frac{5}{4}v^2(\delta_n)\lambda^2(\delta_n)\leq C_2 r_n(x,f).$$
\end{proof}

\noindent The next theorem states that the adaptive estimator as proposed in \eqref{ad_es} mimics the performance of the ideal estimator formed with the locally deterministic bandwidth $\delta_n$.

\begin{theorem} \label{th2}
	For any $f\in B(L)$ with $L>0$ and any point $x\in \mathbb{R}^2$, we have for $n\to \infty$:$$E_f[(f^*(x)-f(x))^2]\leq c(a)v^2(\delta_n)\lambda^2(\delta_n)\leq c^\prime(a)r_n(x,f),$$ where $c(a)$ and $c^{\prime}(a)$ are constants depending on $a$.
\end{theorem}
\begin{proof}
	We decompose the risk in to two parts and consider each part one by one:
	\begin{align*}
	E_f[(f^*(x)-f(x))^2]&=E_f[(f^*(x)-f(x))^2]I(\bar{\delta}\geq \delta_n)\nonumber \\&+E_f[(f^*(x)-f(x))^2]I(\bar{\delta}\leq \delta_n).
	\end{align*}
	\textbf{Case 1: } $\{\bar{\delta}\geq \delta_n\}.$ Note that for any $\delta^{\prime }\geq \delta$ we have $v(\delta)\geq v(\delta^{\prime})$ and $\lambda (\delta)\geq \lambda(\delta^{\prime}).$ Thus it is easy to see that $\psi(\delta^{\prime},\delta)\leq v(\delta)\lambda(\delta)+v(\delta,\delta^{\prime})\lambda(\delta)$. Using the fact that $\int (K_\delta(s)-K_{\delta^{\prime}}(s))^2 ds\leq 2\int K^2_{\delta}(s)ds+2\int K^2_{\delta^{\prime}}(s)ds$, we get that $v^2(\delta,\delta^{\prime})\leq 2[v^2(\delta)+v^2(\delta^{\prime})]\leq 4v^2(\delta)$. Thus, $v(\delta,\delta^{\prime})\leq 2v(\delta)$ which in turn implies, $\psi(\delta,\delta^{\prime})\leq 3v(\delta)\lambda(\delta).$\\
	Now we have a series of inequalities, \begin{align*}
	\lvert f^*(x)-\bar{f}_{\delta_n}(x)\rvert I(\bar{\delta}\geq \delta_n)&\leq \psi(\bar{\delta},\delta_n)\leq \max (\psi(\delta^{\prime},\delta_n):\delta^{\prime}\in \Delta_n,\delta^{\prime}\geq\delta_n)\\&\leq 3v(\delta_n)\lambda(\delta_n).
	\end{align*}
	Thus, 
	\begin{align}\label{case1}
	E_f&[(f^*(x)-f(x))^2]I(\bar{\delta}\geq \delta_n)\nonumber\\
	&=E_f[(f^*(x)-\bar{f}_{\delta_n}(x)+ \bar{f}_{\delta_n}(x)-f_{\delta_n}(x)+f_{\delta_n}(x) - f(x))^2]I(\bar{\delta}\geq \delta_n)\nonumber\\
	&\leq 3[E_f[(f^*(x)-\bar{f}_{\delta_n}(x))^2]I(\bar{\delta}\geq \delta_n)+E_f[(\bar{f}_{\delta_n}(x)-f_{\delta_n}(x))^2]+(f_{\delta_n}(x)-f(x))^2]\nonumber\\
	&\leq 3[9 v^2(\delta_n)\lambda^2(\delta_n)+v^2(\delta_n)+\frac{1}{4}v^2(\delta_n)\lambda^2(\delta_n)] \nonumber \\
	&=c_1v^2(\delta_n)\lambda^2(\delta_n), 
	\end{align}
	where $c_1$ is a constant and where we have used the fact that $\lambda^2(\delta_n)\geq 1.$\\  
	
	\par \noindent \textbf{Case 2: } $\{\bar{\delta}<\delta_n   \}$. Consider the set $B_n(x,\delta,\eta)=\{\lvert \bar{f}_{\delta}(x)-\bar{f}_{\eta}(x)\rvert>\psi(\delta,\eta)   \}$ where $\eta\in \Delta_n, \delta \in \Delta_n$ and $\delta>\eta$. Consider the event $\{\bar{\delta}=\delta/a  \}$ for any $\delta\in \Delta_n$. Since $a>1$, this implies $\delta>\bar{\delta}$. Let ${\Delta_n{(\delta)}}:=\{\eta\in \Delta_n, \eta<\delta\}$. Thus from the definition of $\bar{\delta}$ we get,
	$\{\bar{\delta}=\delta/a \}\displaystyle{\subset \cup_{\eta\in \Delta_n(\delta)}}B_n(x,\delta,\eta).$
	From this it follows, 
	\begin{align*}
	\{\bar{\delta}<\delta_n\}\subset \bigcup \{\bar{\delta}=\delta/a: \delta\in \Delta_n(a\delta_n)  \}\subset \underset{\delta\in \Delta_n(a\delta_n)}{\bigcup}\quad \underset{\eta\in \Delta_n(\delta)}{\bigcup}B_n(x,\delta,\eta).
	\end{align*}
	Thus we have the following series of inequalities,
	\begin{align*}
	E_f\big[(f^*(x)-&f(x))^2I(\bar{\delta}<\delta_n)\big]\leq \sum_{\delta\in \Delta_n(a\delta_n)}E_f[(\bar{f}_{a^{-1}\delta}(x)-f(x))^2I(\bar{\delta}=a^{-1}\delta)]\nonumber\\
	& \quad \quad \quad \quad \quad \quad  \quad \leq \sum_{\delta\in \Delta_n(a\delta_n)} \sum_{\eta\in \Delta_n(\delta)} E_f[(\bar{f}_{a^{-1}\delta}(x)-f(x))^2I(B_n(x,\delta,\eta))].
	\end{align*}
	Since $\delta/a<\delta_n$, $\lvert f_{a^{-1}\delta}(x)-f(x)\rvert\leq v(\delta_n)\lambda(\delta_n)/2\leq v(\delta)\lambda(\delta)/2$. We remark here that the fact that both $\delta$ and $\delta_n$ are in the geometric grid $\Delta_n$ along with the fact that $\delta<a\delta_n$ gives us that $\delta\leq \delta_n$. This explains the rightmost inequality in the above expression. 
	
	\noindent Furthermore, from the definition of $\delta_n$, for any $\eta<\delta\leq\delta_n$ we get,
	\begin{align*}
	\lvert f_{\eta}(x)-f(x)\rvert \leq \frac{v(\delta_n)\lambda(\delta_n)}{2}\leq \frac{v(\delta)\lambda(\delta)}{2}.
	\end{align*} 
	Note that, 
	\begin{align*}
	\lvert \bar{f}_{\delta}(x)-\bar{f}_{\eta}(x)&\rvert =\lvert \bar{f}_{\delta}(x)-\bar{f}_{\eta}(x)-(f_{\delta}(x)-f_{\eta}(x))+f_{\delta}(x)-f(x)+f(x)-f_{\eta}(x) \rvert\nonumber\\
	&\leq \lvert \bar{f}_{\delta}(x)-\bar{f}_{\eta}(x)-(f_{\delta}(x)-f_{\eta}(x)) \rvert+\lvert f_{\delta}(x)-f(x) \rvert+\lvert f(x)-f_{\eta}(x)  \rvert\\
	&\leq \lvert \frac{1}{n}\sum_{i=1}^{n}\zeta_i\rvert +v(\delta)\lambda(\delta),
	\end{align*}
	where $\zeta_i=e^{-\mu x\cdot \theta_i^{\perp}}\bigg(K_{\delta}(\langle x\cdot \theta_i \rangle -s_i )-K_{\delta}(\langle x\cdot \theta_i \rangle -s_i )\bigg)Y_i-(f_{\delta}(x)-f_{\eta}(x))$.
	Thus it follows from the definition of $B_n(x,\delta,\eta)$ and $\psi(\delta,\eta)$that $$B_n(x,\delta,\eta)\subset\Big\{\Big\lvert\frac{1}{n}\sum_{i=1}^{n}\zeta_i\Big\rvert >v(\delta,\eta)\lambda(\eta)  \Big\}.$$
	Therefore, 
	\begin{align}
	&E_f\big[(f^*(x)-f(x))^2I(\bar{\delta}<\delta_n)\big]\nonumber\\ &\leq \sum_{\delta\in \Delta_n(a\delta_n)} \sum_{\eta\in \Delta_n(\delta)} E_f\Big[(\bar{f}_{a^{-1}\delta}(x)-f(x))^2I\big(\{\big\lvert\frac{1}{n}\sum_{i=1}^{n}\zeta_i\big\rvert >v(\delta,\eta)\lambda(\eta)  \}\big)\Big]\nonumber\\
	&\leq \sum_{\delta\in \Delta_n(a\delta_n)} \sum_{\eta\in \Delta_n(\delta)} \bigg(E_f[(\bar{f}_{a^{-1}\delta}(x)-f(x))]^4\bigg)^{\frac{1}{2}}\bigg(P_f\Big(\big\lvert\frac{1}{n}\sum_{i=1}^{n}\zeta_i\big\rvert >v(\delta,\eta)\lambda(\eta) \Big)\bigg)^{\frac{1}{2}}, \nonumber
	\end{align}
	where the last inequality follows on an application of the C-S-B inequality. Also note that
	\begin{align} \label{prev}
	P_f\Big(\Big\lvert\frac{1}{n}\sum_{i=1}^{n}\zeta_i\Big\rvert >v(\delta,\eta)\lambda(\eta) \Big)&\leq P_f\Big(\frac{1}{n}\sum_{i=1}^{n}\zeta_i >v(\delta,\eta)\lambda(\eta)  \Big)\nonumber\\
	& \quad \quad \quad  +P_f\Big(-\frac{1}{n}\sum_{i=1}^{n}\zeta_i >v(\delta,\eta)\lambda(\eta)  \Big). 
	\end{align}
	Let us estimate the first term on the RHS of the previous inequality \eqref{prev}. For this, note that we have {by Markov's inequality}: $$P_f\Big(\frac{1}{n}\sum_{i=1}^{n}\zeta_i>v(\delta,\eta)\lambda(\eta) \Big)\leq E_f\Big(\exp{(\frac{z}{n}\sum_{i=1}^{n}\zeta_i)}\Big)\exp{\Big(-zv(\delta,\eta)\lambda(\eta)\Big)}.$$
	For the i.i.d variables $\zeta_i$ as defined above it is easy to see that $$E_f\Big(\frac{1}{n}\sum_{i=1}^{n}\zeta_i\Big)=0 \text{ and } Var_f(\zeta_i)\leq v^2(\delta,\eta).$$ Furthermore, by {using} the fact that $\zeta_i$ are i.i.d, we can write
	\begin{align}\label{zeta_expect}
	E_f\Big[\exp{\big(\frac{z}{n}\sum_{i=1}^n \zeta_i\big)}\Big]=\Big(E_f\big[\exp{\big(\frac{z}{n} \zeta_1\big)}\big]\Big)^n.
	\end{align}
	We will denote by $K^{\eta}_{\delta}:=e^{-\mu x\cdot \theta^{\perp}}[K_{\delta}(\langle x\cdot \theta \rangle -s )-K_{\eta}(\langle x\cdot \theta \rangle -s )]$.
	In view of equation \eqref{zeta_expect}, {let us} at first evaluate the following conditional expectation:
	\begin{align}
	E_f&[\exp{\big(\frac{z}{n} \zeta_1\big)}|(\theta,s)]\nonumber\\
	&=E_f\Big[\exp{\Big(\frac{z}{n}\big( K^{\delta}_{\eta}(T_\mu f +\epsilon_1)-(f_{\delta}(x)-f_{\eta}(x))\big)\Big)}\big|(\theta,s)\Big] \nonumber\\ 
	&=\exp{\Big(\frac{z}{n}\big( K^{\delta}_{\eta}(T_\mu f(\theta,s))-(f_{\delta}(x)-f_{\eta}(x))\big)\Big)}E_f\Big[\exp{\big(\frac{z}{n}K^{\eta}_{\delta}\epsilon_1\big)}\big|(\theta,s)\Big]\nonumber\\
	&=\exp{\Big(\frac{z}{n}\big(K^{\eta}_{\delta}T_{\mu}f(\theta,s)-(f_{\delta}(x)-f_{\eta}(x))+\frac{z^2\sigma^2}{{2}n^2}(K^{\eta}_{\delta})^2\big)\Big)} \quad (\epsilon_1 \text{ is Gaussian}) \nonumber \\
	&=\exp\Big(\frac{z}{n}(K^{\eta}_{\delta}T_{\mu}f(\theta,s)-(f_{\delta}(x)-f_{\eta}(x))+\frac{z^2\sigma^2}{{2} n^2}{((K^{\eta}_{\delta})^2-E_{(\theta,s)}(K^{\eta}_{\delta})^2)} \nonumber\\
	& \quad \quad \quad \quad \quad \quad \quad \quad \quad +\frac{z^2\sigma^2}{2n^2}E_{(\theta,s)}(K^{\eta}_{\delta})^2\Big)\nonumber\\
	&=\exp{(U_1+U_2)}\exp{\Big(\frac{z^2\sigma^2}{2n^2}E_{(\theta,s)} (K^{\eta}_{\delta})^2\Big)}, \nonumber
	\end{align} 
	where $U_1=\frac{z}{n}(K^{\eta}_{\delta}T_{\mu}f(\theta,s)-(f_{\delta}(x)-f_{\eta}(x))$ and {$U_2=\frac{z^2\sigma^2}{2n^2}((K^{\eta}_{\delta})^2-E_{(\theta,s)}(K^{\eta}_{\delta})^2)$}.
	Observe here that $E_{(\theta,s)}(U_1)=0=E_{(\theta,s)}(U_2)$. Also one can easily verify that $Var_{(\theta,s)}U_1 =Var_{(\theta,s)}(\frac{z}{n}\zeta_1)\leq {(z^2/n) v^2(\delta,\eta)}$. For the calculations below we would also need an estimate on $Var_{(\theta,s)}(U_2)$. Note that $Var_{(\theta,s)}(U_2) =\frac{z^4\sigma^4}{4n^4} Var_{(\theta,s)}(K^{\eta}_{\delta})^2=\frac{z^4\sigma^4}{4n^4} [E_{(\theta,s)}(K^{\eta}_{\delta})^4-(E_{(\theta,s)}(K^{\eta}_{\delta})^2)^2]$. As $\eta<\delta$ and both belong to the geometric grid $\Delta_n$, $E_{\theta,s}(K^{\eta}_{\delta})^2\neq 0$. Thus, $$\frac{E_{(\theta,s)}({K^{\eta}_{\delta}})^4}{(E_{(\theta,s)}(K^{\eta}_{\delta})^2)^2}=\frac{E_{(\theta,s)}({K^{\eta}_{\delta}})^4}{E_{(\theta,s)}({K^{\eta}_{\delta}})^4- Var_{(\theta,s)}(K^{\eta}_{\delta})^2}=\frac{1}{1-\bigg(\frac{Var_{(\theta,s)}(K^{\eta}_{\delta})^2}{E_{(\theta,s)}({K^{\eta}_{\delta}})^4}\bigg)}.$$ Moreover from the fact that $E_{\theta,s}(K^{\eta}_{\delta})^2\neq 0$ it follows that, 
	$$Var_{(\theta,s)}(K^{\eta}_{\delta})^2<E_{(\theta,s)}({K^{\eta}_{\delta}})^4.$$
	Thus we get, $\frac{E_{(\theta,s)}({K^{\eta}_{\delta}})^4}{(E_{(\theta,s)}(K^{\eta}_{\delta})^2)^2}=\tilde{C}(\eta,\delta)>1$ where $\tilde{C}(\eta,\delta)$ is some constant depending upon $\eta$ and $\delta$. This in turn gives us, 
	$$Var_{(\theta,s)}(U_2) \leq\frac{z^4\sigma^4}{4n^4} (\tilde{C}(\eta,\delta)-1)(E_{(\theta,s)}(K^{\eta}_{\delta})^2)^2\leq \frac{z^4\sigma^4}{4n^2} (\tilde{C}(\eta,\delta)-1) v^4(\eta,\delta).$$ Taking $z=\delta \lambda(\eta)/v(\delta,\eta)$ and $\tilde{C}(\eta,\delta)-1=C(\eta,\delta)$ , we get, 
	$$Var_{(\theta,s)}U_2\leq \delta^4\lambda^4(\eta)C(\eta,\delta)\sigma^4/4n^2. $$ 
	Also recall, $E_{(\theta,s)}(K^{\eta}_{\delta})^2\leq \frac{4 \pi n}{4\pi^2L^2+\sigma^2}v^2(\delta,\eta)$. Finally as $U_1$ and $U_2$ are bounded (and thus sub-Gaussian), we get:
	\begin{align*}
	[E_f(\exp{(\frac{z}{n}\zeta_1)})]^n\leq \exp{(\delta^2\lambda^2(\eta))} \exp{\Big(\frac{2\pi\delta^2\lambda^2(\eta)\sigma^2}{4\pi^2L^2+\sigma^2}\Big)}\exp{\Big(\frac{\delta^4\sigma^4\lambda^4(\eta)C(\delta,\eta)}{n}\Big)}.
	\end{align*}
	With $n\to \infty$, we get,
	\begin{align} \label{eq:17}
	[E_f(\exp{(\frac{z}{n}\zeta_1)})]^n\leq \exp{\Big(\delta^2\lambda^2(\eta)\big(1+\big(\frac{{2\pi\sigma^2}}{4\pi^2L^2+\sigma^2}\big)\big) \Big)}. 
	\end{align}
	Using \eqref{eq:17},
	\begin{align}\label{eq:18}
	P_f\Big(\frac{1}{n}\sum_{i=1}^{n}\zeta_i&>v(\delta,\eta)\lambda(\eta) \Big)\leq \exp{\Big(\lambda^2(\eta)\Big(\delta^2\big(1+(\frac{2\pi \sigma^2}{4\pi^2L^2+\sigma^2})\big)-\delta\Big) \Big)}\nonumber \\
	&\leq \exp{(\lambda^2(\eta)(c_1\delta^2-\delta) )} \quad (\text{where $c>1$ is a constant }).
	\end{align}
	Since \eqref{eq:18} is true for all $\delta$, in particular it is true for $\delta=1/2c_1$ and we get, 
	\begin{align*}
	P_f\big(\frac{1}{n}\sum_{i=1}^{n}\zeta_i&>v(\delta,\eta)\lambda(\eta) \big)\leq \exp{\Big(-\frac{\lambda^2(\eta)}{4c_1}\Big)}.
	\end{align*}
	Finally, we get, 
	\begin{align}
	P_f\big(\frac{1}{n}\sum_{i=1}^{n}\zeta_i&>v(\delta,\eta)\lambda(\eta) \big)\leq 2\exp{\Big(-\frac{\lambda^2(\eta)}{4c_1}\Big)}\leq  2\exp{\Big(-\frac{D_2 \log(1/\eta)}{4c_1}\Big)}. \nonumber
	\end{align}
	Now consider,
	\begin{align}
	E_f(\bar{f}_{a^{-1}\delta}(x)-f(x))^4&=E_f(\bar{f}_{a^{-1}\delta}(x)-f_{a^{-1}\delta}(x)+f_{a^{-1}\delta}(x)-f(x))^4\nonumber\\
	&\leq 8E_f[(\bar{f}_{a^{-1}\delta}(x)-f_{a^{-1}\delta}(x))^4]+8(f_{a^{-1}\delta}(x)-f(x))^4\nonumber\\
	&\leq 8E_f[(\frac{1}{n}\sum_{i=1}^{n}z_i)^4]+2v^4(\delta_n)\lambda^4(\delta_n), \nonumber
	\end{align}
	where $z_i=e^{-\mu x\cdot \theta_i^{\perp}} K_{a^{-1}\delta}(\langle \theta_i,x \rangle -s_i)Y_i-f_{a^{-1}\delta}(x)$ are i.i.d. random variables. It is easy to see that $E_f(z_i)=0$ and $Var_f(z_i)\leq nv^2(a^{-1}\delta)$. Thus on expanding $(\frac{1}{n}\sum_{i=1}^{n}z_i)^4$, one can show that,
	\begin{align*}
	8E_f[(\frac{1}{n}\sum_{i=1}^{n}z_i)^4]\leq 8 \bigg[\frac{E_f(z_1)^4}{n^3}+\frac{3 {n \choose 2} n^2 v^4(a^{-1}\delta)}{n^4}\bigg].
	\end{align*}
	As $n\to \infty$, we get,
	\begin{align}
	8E_f[(\frac{1}{n}\sum_{i=1}^{n}z_i)^4]\leq c_3 v^4(a^{-1}\delta), \nonumber
	\end{align}
	where $c_3$ is a positive constant. Recalling that $\eta,\delta\leq \delta_n$, $v^2(\delta/a)=c^*(\delta/a)^{1-2d}/n$ and  $v(\delta_n)\leq v(a^{-1}\delta)$ we have,
	\begin{align}\label{series1}
	E_f&[(f^*(x)-f(x))^2I(\bar{\delta}<\delta_n)] \nonumber \\
	&\leq {c_4} \underset{\delta \in \Delta_n(a\delta_n)}{\sum} \underset{\eta\in \Delta_n(\delta)}{\sum} \bigg(v^4(a^{-1}\delta)+v^4(\delta_n)\lambda^4(\delta_n)\bigg)^{\frac{1}{2}}\exp{\Big(-\frac{\lambda^2(\eta)}{8}\Big)}\nonumber\\
	&\leq c_5  \underset{\delta \in \Delta_n(a\delta_n)}{\sum} \underset{\eta\in \Delta_n(\delta)}{\sum} \frac{\delta^{1-2d}}{n}\lambda^2(\eta)(\eta)^{D_2/8}.
	\end{align}
	The number of elements in the set $\Delta_n$ is less than $\left \lceil{\log n/\log a}\right \rceil =N_n $. Thus, 
	\begin{align}\label{series2}
	\sum_{\eta\in \Delta_n(\delta)} \lambda^2(\eta) \eta^{D_2/8}&\leq \delta^{D_2/8-\alpha}\sum_{\eta\in \Delta_n(\delta)} D_2\log(\frac{1}{\eta})\eta^{\alpha}\nonumber\\
	&\leq \delta^{D_2/8-\alpha} \sum_{j=0}^{N_n}{D_{2}}(\frac{1}{a^{\alpha}})^j\log a,
	\end{align}
	where we use the fact $\eta=a^{-j}$ for some $j$ as $\eta\in \Delta_n$. Since $a>1$, then for a small enough $\alpha$ (the choice to be made precise later) the series in \eqref{series2} converges. Thus $\sum_{\eta\in \Delta_n(\delta)} \lambda^2(\eta) \eta^{D_2/8}\leq c_5(\delta)^{D_2/8-\alpha}$. Finally from \eqref{series1},
	\begin{align}\label{case2}
	E_f[(f^*(x)-f(x))^2I(\bar{\delta}<\delta_n)]&\leq \frac{c_6}{n}\underset{\delta \in \Delta_n(a\delta_n)}{\sum}\delta^{\frac{D_2}{8}-\alpha-2d+1}\nonumber\\
	&\leq \frac{c_6}{n}\sum_{j=0}^{N_n}\bigg(\frac{1}{a^{\alpha}}\bigg)^j\leq c_7 v^2(\delta_n)\lambda^2(\delta_n), 
	\end{align}
	where $\alpha$ is chosen such that $2\alpha\leq D_2/8-2d+1$ and $c_7$ is a constant that depends on $a$.
	The proof of the theorem follows from \eqref{case1} and \eqref{case2}.
\end{proof}
\noindent Now we begin our analysis of the performance of such an adaptive estimator over a class of functions comprising a wide Sobolev regularity scale. At first we show that the rate of convergence in this adaptive procedure is off only by a $\log n$ factor when compared with the minimax optimal rate of estimation for a function $f\in H(\beta)\cap B(L)$. Such a loss of efficiency is in fact ubiquitous in pointwise adaptive estimation of functions (e.g. over Sobolev classes, see \cite{Tsybakov_98}) and can not be done away with. This gives us confidence in the validity of applying the adaptive procedure as proposed in \cite{Cavalier01,Lepski_97} for the adaptive estimation of function from its stochastic ERT data as well.
\begin{theorem}\label{th3}
	Let $L>0$. For any $x\in \mathbb{R}^2$ and for any $\beta>1$ we have:
	$$\limsup_{n\to \infty}\sup_{f\in H(\beta,C)\cap B(L)}\bigg(\frac{n}{\log n}\bigg)^{\frac{2\beta-d}{2\beta+d-1}}E_{f}(f^*(x)-f(x))^2<\infty.$$
\end{theorem}
\begin{proof}
	From equation (9) in \cite{AA_opt_ERT}, we have that for $f\in H(\beta,C)\cap B(L)$, $(f_{\eta}(x)-f(x))^2\leq \tilde{c_8}\eta^{2\beta-2}$ for all $\eta$. Then from the definition $$r_n(x,f)\leq d_1\delta_0^{2\beta-2}c^*\delta_0^{-3}\frac{\log n}{n}$$ where $\delta_n^{-}=\frac{\log n}{n}<\delta_0<1$. If we choose $ \delta_0=(\frac{\log n}{n})^{\frac{1}{2\beta+1}}$ then we get,
	\begin{align}
	r_n(x,f)&\leq d_1\bigg(\frac{\log n}{n}   \bigg)^{\frac{2\beta-2}{2\beta+1}}c^*\bigg(\frac{\log n}{n}\bigg)^{\frac{-3}{2\beta+1}}\bigg(\frac{\log n}{n} \bigg)\nonumber\\
	&=(d_1+c^*)\bigg(\frac{\log n}{n}\bigg)^{\frac{2\beta-2}{2\beta+1}}. \nonumber
	\end{align}
	Now the result follows from Theorem \ref{th2}.
\end{proof}
\begin{remark}
	A function $f\in B(L)$ is said to be locally in $H(\beta,C)$ near $x_0$ if there exists a smooth cut-off function $\phi\in C_c(\mathbb{R}^2)$ with $\phi(x_0)\neq 0$ such that $\phi f \in H(\beta,C)$. With obvious modifications, Theorem \ref{th3} is also true for such functions which are known to be locally in $H(\beta,C)\cap B(L)$ near any arbitrary fixed point $x$.
\end{remark}
\noindent Finally, in the next theorem we show that there can not exist any adaptive estimator on the scale $H(\beta)\cap B(L)$ which is also optimally rate adaptive in the sense of Remark \ref{orc}. Furthermore, we also show that the estimator described by equation \ref{ad_es} achieves the adaptive rate of convergence for $H(\beta)\cap B(L)$. The proof is based on the methods used in \cite{Butucea_00, Butucea_01}.
\begin{theorem} \label{Th4}
	Let $B_n=\{\beta_1<\beta_2<\dots<\beta_{N_n}\}$ such that $\beta_1>1$ is fixed and $\beta_{N_n}\underset{n \to \infty}{\to} \infty$. Then there exists no `optimal adaptive estimator' over the class $H(\beta)\cap B(L)$ for the pointwise risk given by definition \ref{risk_fn}. At the same time the rate $\psi_{n,\beta}=\big(\frac{n}{\log n}\big)^{\frac{\beta-1}{2\beta+1}}$ is the adaptive rate in the sense of definition \ref{adap}.
\end{theorem}
\begin{proof}
	The first part of the theorem essentially claims the non-existence of an adaptive estimator that is optimal across the Sobolev scale. We will establish this result first. Let $\phi_{n,\beta}=n^{\frac{\beta-1}{2\beta+1}}$ denote the optimal minimax rate (up to a constant) of estimation when $\beta$ is known. Suppose that there does exist an adaptive estimator achieving optimal minimax rate for every Sobolev smoothness $\beta$. Then, 
	\begin{align}\label{eq:26}
	C&\geq \limsup_{n\to \infty}\sup_{\beta \in B_n}\sup_{f\in H(\beta,c)\cap B(L)} (\phi_{n,\beta})^{-2}E_{f}(f^*(x)-f(x))^2\nonumber\\
	&\geq\limsup_{n\to \infty}\sup_{\beta \in B_n}\sup_{f\in H(\beta,c)\cap B(L)}\bigg(\frac{\phi_{n,\beta}}{\psi_{n,\beta}}\bigg)^{-2} \psi_{n,\beta}^{-2}E_{f}(f^*(x)-f(x))^2\nonumber\\
	&\geq \liminf_{n\to \infty}\bigg(\sup_{\beta\in B_n}\bigg(\frac{\psi_{n,\beta}}{\phi_{n,\beta}}\bigg)^2\inf_{\hat{f}_n}\sup_{\beta}\sup_{f\in H(\beta,L)}\psi_{n,\beta}^{-2}E_f(\hat{f}_n(x)-f(x))^2\bigg).
	\end{align} We will show that the inequality \eqref{eq:26} gives us a contradiction by showing that the right hand side of the above inequality is unbounded. To that end, we prove the following lemma first.
	\begin{Lemma}\label{lem1}
		For any estimator $\hat{f}_n$ and rate $\psi_{n,\beta}$ as above we have,
		$$\liminf_{n \to \infty}\inf_{\hat{f}_n}\sup_{\beta}\sup_{f\in H(\beta,L)}\psi_{n,\beta}^{-2}E_f(\hat{f}_n(x)-f(x))^2\geq 1.$$
	\end{Lemma}
	\begin{proof}
		We will denote $\beta_{N_n}$ by $\beta_{N}$. Consider the two hypotheses, $f_{n,0}(x)=0$ and $f_{n,1}(x)=Ah^{\beta_1-1}\eta((x-x_0)/h)$ where $0<A<1$ is a constant, $h=(\frac{\log n}{n})^{\frac{1}{2\beta_1+1}}$ where $\eta(x)$ is a compactly supported function in $H(\beta,L)$ such that $\eta(0)=1$. It is easy to show that $f_{n,1}(x)\in H(\beta_1,L)$ and $f_{n,0}\in H(\beta_N,L)$. Furthermore, one can easily find $\delta$ such that $d(f_{n,1},f_{n,0})>1-2\delta$ for some $0<\delta<1/2$. Note the inequality,
		\begin{align}\label{ineqN}
		\inf_{\hat{f}_n}\sup_{\beta}\sup_{f\in H(\beta,L)}\psi_{n,\beta}^{-2}E_f(\hat{f}_n(x)-f(x))^2\geq\inf_{\hat{f}_n}\max &\bigg\{E_{f_{n,0}}[\psi_{n,\beta_{N}}^{-2} \lvert\hat{f}_n(x_0)-f_{n,0}(x_0)\rvert^2],\nonumber\\
		&E_{f_{n,1}}[\psi_{n,\beta_{1}}^{-2} \lvert\hat{f}_n(x_0)-f_{n,1}(x_0)\rvert^2]\bigg\}.
		\end{align}
		To establish a lower bound for the expression on the right hand side of the inequality \eqref{ineqN} above, we will use \cite{Tsybakov98}[A1, Theorem 6]. Let $E_0[\cdot]:=E_{f_{n,0}}[\cdot]$ and $E_1[\cdot]:=E_{f_{n,1}}[\cdot]$ and the associated probability measures be denoted by $P_0$ and $P_1$ respectively. To apply \cite{Tsybakov98}[A1, Theorem 6], we will need to show that for $0<\alpha<1$ and $\tau>0$, $P_1[\frac{dP_0}{dP_1}\geq \tau]\geq1-\alpha$. 
		\begin{align}
		P_1\bigg[\frac{dP_0}{dP_1}\geq \tau\bigg]&=P_1\bigg[\prod_{i=1}^{n}\frac{p_0(Y_i)}{p_1(Y_i)}\geq \tau\bigg]\nonumber\\
		&=P_1\bigg[\frac{1}{\sqrt{\log n}}\sum_{i=1}^{n}\log\frac{p_0(Y_i)}{p_1(Y_i)}\geq \frac{\log\tau}{\sqrt{\log n}}\bigg].
		\end{align}
		Let $Z_{n,i}=\frac{1}{\sqrt{\log n}}\log\frac{p_0(Y_i)}{p_1(Y_i)}$ which are i.i.d. random variables. Let $E_1[Z_{n,i}]$ and $V_1[Z_{n,i}]$ be the expectation and variance $Z_{n,i}$ with respect to the probability measure $P_1$ corresponding to the experiments with observations generated by $f_{n,1}$. Let $\sigma_n=\sum_{i=1}^n V_1[Z_{n,i}]$ and define $U_{n,i}=(Z_{n,i}-E_1[Z_{n,i}])/\sigma_n$.
		We state the following lemma, the proof for which will be given in the appendix.
		
		\begin{Lemma} \label{aux_lemma}
			For $Z_{n,i}$ and $U_{n,i}$ as defined above, we have,
			\begin{enumerate}
				\item[(a)] $\sum_{i=1}^{n}E_1[Z_{n,i}]\geq -c_8\sqrt{\log n}$ where $c_8>0$ is a small enough constant.
				\item[(b)] $\sigma_n^2\geq c_9>0$.
				\item[(c)] $\lim_{n\to\infty}\sum_{i=1}^n E_1[\lvert U_{n,i}\rvert^3]=0$.
			\end{enumerate}
		\end{Lemma}
		\noindent From part (c) in Lemma \ref{aux_lemma}, we conclude that $U_n=\sum_{i=1}^{n}U_{n,i}$ converges in law to the Normal Distribution $N(0,1)$ (Lyapunov's CLT).
		Thus we can rewrite $P_1[\sum_{i=1}^nZ_{n,i}\geq \frac{\log \tau}{\sqrt{\log n}}]=P_1[U_n\geq m_n]$ where $m_n=\frac{\frac{\log \tau}{\sqrt{\log n}}-\sum_{i=1}^nE_1[Z_{n,i}]}{\sigma_n}$. Now choose $\tau=n^{-r}$ for some $r>c_8>0$ (the choice of $r$ will be made precise below), then using (a,b) of Lemma \ref{aux_lemma}, we get $m_n\leq \frac{-r\sqrt{\log n}+c_8\sqrt{\log n}}{\sqrt{c_9}}\to - \infty$ as $n\to \infty$. This shows that, $$P_1\bigg[\frac{dP_0}{dP_1}\geq \tau\bigg]=P_1\bigg[\sum_{i=1}^nZ_{n,i}\geq \frac{\log \tau}{\sqrt{\log n}}\bigg]=P_1[U_n\geq m_n]\to 1 \text{ as } n \to \infty.$$
		
		Furthermore, from \cite[A1 Theorem 6]{Tsybakov_98}, it follows that for $q_n>0,\tau>0, 0<\alpha<1$ and $0<\delta<1/2$:
		\begin{align*}
		\inf_{\hat{f}_n}\max &\bigg\{E_{f_{n,0}}[\psi_{n,\beta_{N}}^{-2} \lvert\hat{f}_n(x_0)-f_{n,0}(x_0)\rvert^2]
		,E_{f_{n,1}}[\psi_{n,\beta_{1}}^{-2} \lvert\hat{f}_n(x_0)-f_{n,1}(x_0)\rvert^2]\bigg\}\\
		&\geq  \frac{(1-\alpha)(1-2\delta)^2\tau q_n^2\delta^2}{(1-2\delta)^2+\tau q_n^2\delta^2}.
		\end{align*} Thus, to satisfy the required lower bound, take $q_n=\frac{\psi_{n,\beta_1}}{\psi_{n,\beta_{N}}}$. Then 
		\begin{align}\label{eq:29}
		\liminf_{n \to \infty} \tau q_n^2&=\liminf_{n\to\infty}n^{-r}\bigg(\frac{\log n}{n}\bigg)^{\frac{2\beta_1-2}{2\beta_1+1}-\frac{2\beta_N-2}{2\beta_N+1}}\nonumber\\
		&=\liminf_{n\to\infty}n^{-r}\bigg(\frac{\log n}{n}\bigg)^{\frac{6(\beta_1-\beta_N)}{(2\beta_1+1 )(2\beta_N+1)}}.
		\end{align}
		Thus if we choose $r$ such that $\frac{6(\beta_N-\beta_1)}{(2\beta_1+1 )(2\beta_N+1)}>r>c_8>0$, then $$\liminf_{n \to \infty} \tau q_n^2 \to \infty \text{ as }n\to \infty.$$
		Consider then,
		\begin{align}
		\liminf_{n \to \infty} \frac{(1-\alpha)(1-2\delta)^2\tau q_n^2\delta^2}{(1-2\delta)^2+\tau q_n^2\delta^2} \to 1,
		\end{align}
		as $\delta$ and $\alpha$ can be chosen as small as desired. Thus from \cite{Tsybakov98} [A1, Theorem 6], we get:
		\begin{align*}
		\inf_{\hat{f}_n}\max &\bigg\{E_{f_{n,0}}[\psi_{n,\beta_{N}}^{-2} \lvert\hat{f}_n(x_0)-f_{n,0}(x_0)\rvert^2]
		,E_{f_{n,1}}[\psi_{n,\beta_{1}}^{-2} \lvert\hat{f}_n(x_0)-f_{n,1}(x_0)\rvert^2]\bigg\}\geq 1.
		\end{align*}
		This concludes the proof of Lemma \ref{lem1}.
	\end{proof}
	\noindent Coming back to the proof of Theorem \ref{Th4}, first observe that $$\sup_{\beta\in B_n}\bigg(\frac{\psi_{n,\beta}}{\phi_{n,\beta}}\bigg)^2 \underset{n \to \infty}{\to} \infty.$$ This along with lemma \ref{lem1} and inequality \eqref{eq:26} gives us a contradiction showing that there can not exist an adaptive estimator that achieves the optimal minimax rate for all $\beta\in B_n.$
	\par \noindent Now we prove the second part of Theorem \ref{Th4}. Let there be another adaptive estimator $f^{**}(x)$ and another sequence $\gamma_{n,\beta}$ such that for any $x\in \mathbb{R}^2$ and for any $\beta>1$ we have,
	$$\limsup_{n\to \infty}\sup_{f\in H(\beta,c)\cap B(L)}(\gamma_{n,\beta})^{-2}E_{f}(f^{**}(x)-f(x))^2<\infty.$$
	\noindent Furthermore let there exist $ \beta^{\prime}$ such that $\frac{\gamma_{n,\beta^{\prime}}}{\psi_{n,\beta^{\prime}}}\underset{n \to \infty}{\to} 0.$
	First of all, {$\gamma_{n,\beta^{\prime}}\geq (\frac{1}{n})^{\frac{\beta^{\prime}-1}{2\beta^{\prime}+1}}.$} Now define: $\kappa_n^{r^\prime}=(\frac{1}{n})^{\frac{\beta^{\prime}-(1+\epsilon_0)}{2\beta^{\prime}+1}}$ where $\epsilon_0>$ is small. Take any $\beta^{\prime\prime}>\beta^{\prime}$. We claim that, 
	\begin{align} \label{claim1}
	\liminf_{n \to \infty} \frac{\gamma_{n,\beta^{\prime\prime}}}{\kappa_n^{r^\prime}}=\infty.
	\end{align}
	{We will prove this claim later, but for the moment if we assume this claim is true and consider,}
	\begin{align}
	\frac{\gamma_{n,\beta^{\prime}}}{\psi_{n,\beta^{\prime}}} \frac{\gamma_{n,\beta^{\prime\prime}}}{\psi_{n,\beta^{\prime\prime}}}&\geq \bigg(\frac{1}{\log n}\bigg)^{\frac{\beta^{\prime}-1}{2\beta^{\prime}+1}}\frac{\gamma_{n,\beta^{\prime\prime}}}{\kappa_n^{r^\prime}}\frac{\kappa_n^{r^\prime}}{\bigg(\frac{\log n}{n}\bigg)^{\frac{\beta^{\prime\prime}-1}{2\beta^{\prime\prime}+1}}}\nonumber\\
	&\geq \frac{\gamma_{n,\beta^{\prime\prime}}}{\kappa_n^{r^\prime}} n^{\frac{\beta^{\prime\prime}-1}{2\beta^{\prime\prime}+1}-\frac{\beta^{\prime}-(1+\epsilon_0)}{2\beta^{\prime}+1}}\bigg(\frac{1}{\log n}\bigg)^{\frac{\beta^{\prime}-1}{2\beta^{\prime}+1}+\frac{\beta^{\prime\prime}-1}{2\beta^{\prime\prime}+1}} \to \infty.
	\end{align}
	This follows from the assumption in the claim made above (see \eqref{claim1}) and the fact that $\frac{n^{\alpha_1}}{(\log n)^{\alpha_2}}\to \infty$ as long as $\alpha_1>0$ and $\alpha_2>0$. (Note also that ${\frac{\beta^{\prime\prime}-1}{2\beta^{\prime\prime}+1}-\frac{\beta^{\prime}-(1+\epsilon_0)}{2\beta^{\prime}+1}}>0.$) Thus the only thing that remains to show is that (\ref{claim1}) holds. To that end, assume $\liminf_{n \to \infty} \frac{\gamma_{n,\beta^{\prime\prime}}}{\kappa_n^{r^\prime}}\leq C<\infty.$ We will show that this gives rise to a contradiction. Recall,
	\begin{align}\label{eq:25}
	&\limsup_{n\to \infty}\sup_{\beta\in B_n}\sup_{f\in H(\beta,c)\cap B(L)}(\gamma_{n,\beta})^{-2}E_{f}(f^{**}(x)-f(x))^2\leq C^{*}<\infty\nonumber\\
	\implies &\limsup_{n\to \infty} \sup_{f\in H(\beta,c)\cap B(L)} \max\bigg\{ \bigg(\frac{\psi_{n,\beta^{\prime}}}{\gamma_{n,\beta^{\prime}}}\bigg)^2 (\psi_{n,\beta})^{-2}E_{f}(f^{**}(x)-f(x))^2,\nonumber\\ 
	&\bigg(\frac{\kappa_{n}^{r^\prime-r}\kappa_{n}^{r}}{\gamma_{n,\beta^{\prime\prime}}}\bigg)^2 (\kappa_{n}^{r^\prime})^{-2}E_{f}(f^{**}(x)-f(x))^2\bigg\}\leq C^* \quad \quad (r^{\prime}>r) \nonumber\\
	\implies & \limsup_{n\to \infty} \min\biggl\{\bigg(\frac{\psi_{n,\beta^{\prime}}}{\gamma_{n,\beta^{\prime}}}\bigg)^2,\bigg(\frac{\kappa_{n}^{r^\prime-r}\kappa_{n}^{r}}{\gamma_{n,\beta^{\prime\prime}}}\bigg)^2 \biggr\} \cdot  \liminf_{n \to \infty} \inf_{\hat{f}} \sup_{f\in H(\beta,c)\cap B(L)} \nonumber\\&\max\biggl\{(\psi_{n,\beta})^{-2}E_{f}(\hat{f}(x)-f(x))^2,
	(\kappa_{n}^{r^\prime})^{-2}E_{f}(\hat{f}(x)-f(x))^2\biggr\}\leq C^*.
	\end{align}
	Since $r^{\prime}>r$, we have $\bigg(\frac{\kappa_{n}^{r^\prime-r}\kappa_{n}^{r}}{\gamma_{n,\beta^{\prime\prime}}}\bigg)\to \infty$ and $\bigg(\frac{\psi_{n,\beta^{\prime}}}{\gamma_{n,\beta^{\prime}}}\bigg)\to \infty$ by hypothesis. Finally, 
	$$\liminf_{n \to \infty} \inf_{\hat{f}} \sup_{f\in H(\beta,c)\cap B(L)} \max\biggl\{\frac{E_{1}(\hat{f}(x)-f_{n,1}(x))^2}{(\psi_{n,\beta})^{2}},
	\frac{E_{0}(\hat{f}(x)-f_{n,0}(x))^2}{(\kappa_{n}^{r^\prime})^{2}}\biggr\}\geq 1,$$ similar to what was done while proving a lower bound for inequality \eqref{ineqN}. The only change is that we consider now $q_n=\frac{\psi_{n,\beta^{\prime}}}{\kappa_{n}^{r^\prime}}$ while proving the relation \eqref{eq:29}. This gives us a contradiction by showing that the left-hand side of inequality \eqref{eq:25} is $\infty$. 
\end{proof}
\section{Appendix}
\begin{proof}[Proof of Lemma \ref{aux_lemma}]
	\textbf{(a)} Let the distribution function for noise be given by $p_{\epsilon}(u)=\frac{1}{\sqrt{2\pi\sigma^2}}e^{\frac{-u^2}{2\sigma^2}}$ . For the proof of this part, first consider, \begin{align}
	E_1[&Z_{n,i}]=\frac{1}{\sqrt{\log n}}E_{(\theta,s)}\bigg[E_{1|(\theta,s)}\bigg[\log \frac{p_{\epsilon}(Y_i)}{p_{\epsilon}(Y_i-T_{\mu}f_{n,1}(\theta_i,s_i))}\bigg]\bigg]\nonumber\\
	&=\frac{-1}{\sqrt{\log n}}E_{(\theta,s)}\bigg[\int \log \frac{p_{\epsilon}(Y_i-T_{\mu}f_{n,1}(\theta_i,s_i))}{p_{\epsilon}(Y_i)} p_{\epsilon}(Y_i-T_{\mu}f_{n,1}(\theta_i,s_i)) dY_i\bigg]\nonumber\\
	&\geq \frac{-1}{\sqrt{\log n}}E_{\theta,s}(T_{\mu}f_{n,1}(\theta_i,s_1))^2.  \nonumber
	\end{align} 
	Recall that for $f_{n,1}=Ah^{\beta_1-1}\eta((x-x_0)/h)$ where $h=\bigg(\frac{\log n}{n}\bigg)^{\frac{1}{2\beta+1}}$, similar  to equation (18) in \cite{AA_opt_ERT} we have, $\int_{Z}(T_{\mu}f_{n,1}(\theta_i,s_i))^2 ds d\theta\leq c_8 h^{2\beta+1}$ where $c_8$ is a constant that can be made as small as desired by choosing a small enough $A$. In particular, we will choose $A$ such that $\frac{6(\beta_N-\beta_1)}{(2\beta_1+1 )(2\beta_N+1)}>c_8>0$. We remark here that in deriving the estimate for $\int_{Z}(T_{\mu}f_{n,1}(\theta_i,s_i))^2 ds d\theta$ as above, we assume that the design points satisfy a certain feasibility condition (\cite[Assumption B2]{AA_opt_ERT}): $E_{(\theta,s)}\bigg[\sum\limits_{i=1}^{n}g(\theta_i,s_i)\bigg]\leq C_3\int\limits_{Z}g(\theta,s)dsd\theta.$
	Thus
	$$\sum_{i=1}^{n}E_1[Z_{n,i}]\geq \frac{-1}{\sqrt{\log n}}n E_{\theta,s}(T_{\mu}f_{n,1}(\theta_i,s_1))^2  \geq-c_8\sqrt{\log n}.$$ 
	\par\noindent \textbf{Proof of part (b)} We want to show that $\sigma_n^2=\sum_{i=1}^{n}V_1[Z_{n,i}]$ is bounded below. First note that from the `law of total variance' $V_{1}[Z_{n,i}]\geq E_{(\theta,s)}[V_{1|(\theta,s)}[Z_{n,i}]]$. Consider
	\begin{align}
	Var_{1|(\theta,s)}[Z_{n,i}]&=\frac{1}{\log n}Var_{1|(\theta,s)}\bigg[\log \frac{p_{\epsilon}(Y_i)}{p_{\epsilon}(Y_i-T_{\mu}f_{n,1}(\theta_i,s_i))}\bigg]\nonumber\\
	&=\frac{1}{\log n}\bigg[E_{1|(\theta,s)}\bigg[\log^2 \frac{p_{\epsilon}(Y_i)}{p_{\epsilon}(Y_i-T_{\mu}f_{n,1}(\theta_i,s_i))}\bigg]\nonumber\\
	&- \bigg(E_{1|(\theta,s)}\bigg[\log \frac{p_{\epsilon}(Y_i)}{p_{\epsilon}(Y_i-T_{\mu}f_{n,1}(\theta_i,s_i))}\bigg] \bigg)^2\bigg]. \nonumber
	\end{align}
	Recall that noise has been assumed to have a Gaussian distribution $\sim N(0,\sigma^2)$  Thus,
	\begin{align}\label{eq:30}
	&E_{1|(\theta,s)}\bigg[\log^2 \frac{p_{\epsilon}(Y_i)}{p_{\epsilon}(Y_i-T_{\mu}f_{n,1}(\theta_i,s_i))}\bigg]\nonumber\\&=\frac{1}{\sqrt{2\pi \sigma^2}}\int\log^2 \exp{\Big(\frac{(Y_i-T_{\mu}f_{n,1}(\theta_i,s_i))^2-Y_i^2}{2\sigma^2}\Big)}\exp{\Big(\frac{(Y_i-T_{\mu}f_{n,1}(\theta_i,s_i))^2}{2\sigma^2}\Big)}dY_i\nonumber\\
	&=\frac{1}{4\sigma^4\sqrt{2\pi\sigma^2}}\int \bigg(T_{\mu}^4(f_{n,1}(\theta_i,s_i)+4Y_i^2T_{\mu}^2f_{n,1}(\theta_i,s_i))-4Y_iT_{\mu}^3f_{n,1}(\theta_i,s_i)\bigg)\nonumber\\
	&\quad \quad \quad \quad \quad \quad \quad \exp{\Big(\frac{(Y_i-T_{\mu}f_{n,1}(\theta_i,s_i))^2}{2\sigma^2}\Big)}dY_i\nonumber\\
	&=\frac{1}{{4}\sigma^4}\bigg(T_{\mu}^4f_{n,1}(\theta_i,s_i)+ {4}\sigma^2T_{\mu}^2f_{n,1}(\theta_i,s_i)\bigg). 
	\end{align}
	On the other hand,
	\begin{align}
	\bigg(E_{1|(\theta,s)} \bigg[\log \frac{p_{\epsilon}(Y_i)}{p_{\epsilon}(Y_i-T_{\mu}f_{n,1}(\theta_i,s_i))}\bigg]\bigg)^2
	=\frac{T_{\mu}^4f_{n,1}(\theta_i,s_i)}{4\sigma^4}. \nonumber
	\end{align}
	Thus $Var_{1|(\theta,s)}[Z_{n,i}]=\frac{{4}(T_{\mu}f_{n,1}(\theta_i,s_i))^2}{\log n\sigma^2}$ and hence, $$\sum_{i=1}^{n}E_{\theta,s}[Var_{1|(\theta,s)}[Z_{n,i}]]=\frac{n}{\sigma^2\log n}\int_{Z}(T_{\mu}f_{n,1}(\theta_i,s_i))^2dsd\theta=c_{10}\frac{n}{\log n}h^{2\beta+1}>0.$$
	\par 
	
	\textbf{Proof of part (c)}\begin{align}\label{eq:32}
	E_1\lvert U_{n,i}^3\rvert& =\frac{1}{\sigma_n^3}E_1\lvert Z_{n,i}^3-(E_1[Z_{n,i}])^3 -3 (Z_{n,i})^2E_1[Z_{n,i}]+3(Z_{n,i})(E_1[Z_{n,i}])^2  \rvert\nonumber\\
	&\leq \frac{1}{\sigma_n^3}[E_1\lvert Z_{n,i}\rvert^3+(E_1\lvert Z_{n,i}\rvert)^3 +3 E_1\lvert Z_{n,i} \rvert^2E_1\lvert Z_{n,i} \rvert+3E_1\lvert Z_{n,i} \rvert(E_1\lvert Z_{n,i} \rvert)^2 ]\nonumber \\
	&\leq \frac{1}{\sigma_n^3}[E_1\lvert Z_{n,i}\rvert^3+(E_1\lvert Z_{n,i} \rvert)^3 +3 E_1\lvert Z_{n,i} \rvert^2E_1\lvert Z_{n,i} \rvert+(E_1\lvert Z_{n,i} \rvert)^3].  
	\end{align}
	Now we consider each of the above terms one by one. First of all $E_1\lvert Z_{n,i}\rvert=E_{\theta,s}[E_{1|(\theta,s)}\lvert Z_{n,i}\rvert] $. Thus using Pinsker's second inequality to calculate:
	\begin{align}\label{eq:33}
	E_{1|(\theta,s)}\lvert Z_{n,i}\rvert &=\frac{1}{\sqrt{\log n}}\int\bigg\lvert \log \frac{p_{\epsilon}(Y_i)}{p_{\epsilon}(Y_i-T_{\mu}f_{n,1}(\theta_i,s_i))}\bigg\rvert p_{\epsilon}(Y_i-T_{\mu}f_{n,1}(\theta_i,s_i)) dY_i \nonumber \\
	&\leq \frac{1}{\sqrt{\log n}}\bigg[\frac{T_{\mu}f_{n,1}(\theta_i,s_i)}{\sigma}+\frac{T_{\mu}^2f_{n,1}(\theta_i,s_i)}{2\sigma^2}\bigg] \cite{Tsybakov_book}\text{[Lemma 2.5]}. 
	\end{align}
	Also note that since the cylinder $Z=[0,2\pi]\times[-1,1]$ has finite measure, we have:
	\begin{align}\label{eq:34}
	\lvert\int_{Z} T_{\mu}f_{n,1}(\theta,s) dsd\theta\rvert\leq \int_Z \lvert T_{\mu}f_{n,1}(\theta,s)\rvert dsd\theta\leq c_{10}\bigg(\int_Z\lvert T_{\mu}f_{n,1}(\theta,s)\rvert^2 dsd\theta\bigg)^{1/2}. 
	\end{align}
	From inequalities \eqref{eq:33} and \eqref{eq:34}, we get:
	\begin{align}
	E_1\rvert{Z_{n,i}}\lvert &\leq \frac{c_{11}}{\sqrt{\log n}}\bigg[\bigg(\frac{\log n}{n}\bigg)^{1/2}+\bigg(\frac{\log n}{n}\bigg)\bigg]\nonumber\\
	&\leq \frac{c_{12}}{\sqrt{\log n}}\bigg(\frac{\log n}{n}\bigg)^{1/2} \quad \quad (0<\frac{\log n}{n}\leq \bigg(\frac{\log n}{n}\bigg)^{1/2}<1 \text{ for } n\geq 3).\nonumber
	\end{align}
	Finally,
	\begin{align}\label{eq:35}
	\sum_{i=1}^{n}\bigg(E_1\rvert{Z_{n,i}}\lvert\bigg)^3\leq c_{12}n\bigg(\frac{1}{n}\bigg)^{3/2}\to 0 \quad \text{as }n\to\infty. 
	\end{align}
	Using \eqref{eq:30}, we have:
	$$E_{1|(\theta,s)}[\rvert{Z_{n,i}}\lvert^2]\leq \frac{1}{{4}\sigma^4\log n}\bigg(T_{\mu}^4f_{n,1}(\theta_i,s_i)+ {4}\sigma^2T_{\mu}^2f_{n,1}(\theta_i,s_i)\bigg).$$
	
	Then using the fact that {$\bigg\lvert T_{\mu}f_{n,1}(\theta_i,s_i)\bigg\rvert \leq c_{13}h^{{\beta}}= c_{13}\bigg(\frac{\log n}{n}\bigg)^{\frac{\beta}{2\beta+1}}$},
	\begin{align}
	&E_{1}[\rvert{Z_{n,i}}\lvert^2]\leq \frac{c_{14}}{\log n}\bigg[\bigg(\frac{\log n}{n}\bigg)^{\frac{4\beta}{2\beta+1}}+\bigg(\frac{\log n}{n}\bigg)^{\frac{2\beta}{2\beta+1}}\bigg]\nonumber\\
	&\leq \frac{c_{15}}{\log n}\bigg[\bigg(\frac{\log n}{n}\bigg)^{\frac{2\beta}{2\beta+1}}\bigg]  \quad \quad (0<\bigg(\frac{\log n}{n}\bigg)^{\frac{4\beta}{2\beta+1}}\leq \bigg(\frac{\log n}{n}\bigg)^{\frac{2\beta}{2\beta+1}}<1 \text{ for } n\geq 3). \nonumber
	\end{align}
	Finally,
	\begin{align}\label{eq:37}
	\sum_{i=1}^{n}E_{1}\lvert Z_{n,i}\rvert^2 E_{1}\lvert Z_{n,i}\rvert &\leq c_{16}\frac{n}{\log n}\frac{1}{\sqrt{\log n}}\bigg(\frac{\log n}{n}\bigg)^{\frac{6\beta+1}{4\beta+2}}\nonumber\\
	&\leq c_{16}\frac{1}{\sqrt{\log n}}\bigg(\frac{\log n}{n}\bigg)^{\frac{2\beta-1}{4\beta+2}} \quad \to 0 \text{ as }n\to \infty. 
	\end{align}
	Now we consider $\sum_{i=1}^n E_1\lvert Z_{n,i}\rvert^3$. For this we first evaluate:
	\begin{align}
	&E_{1|(\theta,s)}\lvert Z_{n,i}\rvert^3=\frac{1}{(\log n)^{3/2}}\int \bigg\lvert \log\frac{p_\xi(Y_i)}{p_{\epsilon}(Y_i-T_{\mu}f_{n,1}(\theta_i,s_i))} \bigg\rvert^3 p_{\epsilon}(Y_i-T_{\mu}f_{n,1}(\theta_i,s_i)) dY_i\nonumber\\
	&\leq \frac{1}{(\log n)^{3/2}}\int\bigg[\bigg\lvert  T_{\mu}^2f_{n,1}(\theta_i,s_i)-2Y_iT_{\mu}f_{n,1}(\theta_i,s_i) \bigg\rvert  \bigg]^3p_{\epsilon}(Y_i-T_{\mu}f_{n,1}(\theta_i,s_i)) dY_i\nonumber\\
	&\leq \frac{1}{(\log n)^{3/2}} \int\bigg[\bigg\lvert T^6_{\mu}f_{n,1}(\theta_i,s_i) \bigg\rvert+8 \bigg\lvert Y_i \bigg\rvert^3\bigg\lvert T_{\mu}^3f_{n,1}(\theta_i,s_i)\bigg\rvert+12\bigg\lvert Y_i \bigg\rvert^2\bigg\lvert T_{\mu}^4f_{n,1}(\theta_i,s_i) \bigg\rvert\nonumber\\+&6\bigg\lvert Y_i\bigg\rvert\bigg\lvert T_{\mu}^5f_{n,1}(\theta_i,s_i) \bigg\rvert \bigg]\frac{\exp{\big(-(Y_i-T_{\mu}f_{n,1}(\theta_i,s_i))^2/2\sigma^2}\big)}{\sqrt{2\pi\sigma^2}} dY_i\nonumber\\
	& \leq \frac{c_{16}}{(\log n)^{3/2}}\bigg(\frac{\log n}{n}\bigg)^{\frac{6\beta}{2\beta+1}}, \nonumber
	\end{align}
	where the last inequality follows from the previous one by integrating each term and using the fact that $\lvert T_{\mu}f_{n,1}(\theta_i,s_i)\rvert\leq c_{13}\bigg(\frac{\log n}{n}\bigg)^{\frac{\beta}{2\beta+1}}$. Thus:
	\begin{align}\label{eq:39}
	\sum_{i=1}^n E_1 \lvert Z_{n,i}\rvert^3&\leq \frac{c_{17}}{(\log n)^{1/2}} \frac{ n}{\log n}\bigg(\frac{\log n}{n}\bigg)^{\frac{6\beta}{2\beta+1}}\nonumber\\
	&\leq  \frac{c_{17}}{(\log n)^{1/2}}\bigg(\frac{\log n}{n}\bigg)^{\frac{4\beta-1}{2\beta+1}} \quad \to 0 \text{ as }n\to \infty. 
	\end{align}
	Equations \ref{eq:32},\ref{eq:35},\ref{eq:37} and \ref{eq:39} together prove part (c).
\end{proof}

\bibliographystyle{plain}
\bibliography{refs_2_stat}

\end{document}